\documentclass[11pt]{article}

\usepackage[cp1252]{inputenc} 
\usepackage[a4paper]{geometry}
\usepackage{lmodern}
\usepackage[T1]{fontenc}
\usepackage{amsmath,amssymb,amsthm}
\usepackage{graphicx}
\usepackage{textcomp}
\usepackage{enumerate}

\usepackage[usenames,dvipsnames]{pstricks}
\usepackage{epsfig}
\usepackage{pst-grad} 
\usepackage{pst-plot} 

\newtheorem{prop}{Proposition}[section]

\newtheorem{lem}[prop]{Lemma}

\newtheorem{thm}{Theorem}[section]

\theoremstyle{definition}
\newtheorem{defi}[prop]{Definition}

\newtheorem{rk}[prop]{Remark}

\newcommand{\R}{\mathbb{R}}
\newcommand{\Q}{\mathbb{Q}}
\newcommand{\Z}{\mathbb{Z}}
\newcommand{\N}{\mathbb{N}}

\font\dsrom=dsrom10 scaled 1200
\def \indic{\textrm{\dsrom{1}}}

\newcommand{\Addresses}{{
  \bigskip
  \footnotesize

\textit{E-mail address}, M.~Bouljihad: \texttt{mohamed.bouljihad@ens-lyon.fr\\}
\textsc{UMPA, ENS Lyon, 46 all\'ee d'Italie, 69007 Lyon, France
   \\ IRMAR, 263 avenue du G\'en\'eral Leclerc, 35000 Rennes, France}\par\nopagebreak

}}

\title{Existence of rigid actions for finitely-generated non-amenable linear groups}

\author{Mohamed Bouljihad}

\date{}

\begin{document}
\maketitle

\begin{abstract}
We show that every finitely-generated non-amenable linear group over a field of characteristic zero admits an ergodic action which is rigid in the sense of Popa. If this group has trivial solvable radical, we prove that these actions can be chosen to be free. Moreover, we give a positive answer to a question raised by Ioana and Shalom concerning the existence of such a free action for $\mathbb{F}_2\times\Z$. More generally, we show that for groups $\Gamma$ considered by Fernos in \cite{Fer}, e.g. Zariski-dense subgroups in PSL$_n(\R)$, the product group $\Gamma\times \Z$ admit a free ergodic rigid action. 
Also, we investigate how rigidity of an action behaves under restriction and co-induction. In particular, we show that rigidity passes to co-amenable subgroups.
\end{abstract}


\setlength{\parindent}{0cm}

\section{Introduction}

In his famous three pages paper \cite{Kaz}, Kazhdan introduced the property (T) for locally compact groups to prove that certain lattices are finitely generated. This notion was first defined in terms of representations and has since become a basic notion in various areas of mathematics such as group theory, ergodic theory, and operator algebras.

In this very same paper, Kazhdan implicitely introduced the notion of relative property (T) for pair of groups $(G,H)$ where $H\subset G$, showing that the pair $(\R^2\rtimes$SL$_2(\R),\R^2)$ has the relative property (T). Yet, it is Margulis \cite{Mar} who really coined this term when he constructed the first explicit examples of families of expander graphs using the relative property (T) for the pair $(\Z^2\rtimes$SL$_2(\Z),\Z^2)$.

Lately, in the spirit of relative property (T) for pair of groups, Popa \cite{Pop} introduced the notion of relative property (T) for inclusions $B\subset N$ of von Neumann algebras (also called rigidity of inclusion).

Given a free ergodic action of a countable group $\Gamma$ on an atomless standard probability space $(X,\mu)$ by measure preserving transformations, there is an associated canonical von Neumann algebra; the group measure space construction L$^\infty(X)\rtimes\Gamma$ (\cite{MvN}). Following \cite{Gab}, we say that the action $\Gamma\curvearrowright (X,\mu)$ has \emph{the property (T) relative to the space} if the inclusion L$^\infty(X)\subset$ L$^\infty(X)\rtimes\Gamma$ has the relative property (T). This notion first appeared in \cite{Pop}, where such actions are called rigid. For an ergodic theoretical equivalent formulation of this property, see Proposition \ref{propIoa} below.

In the case of group actions by automorphims on abelian groups, this notion is linked to the Kazhdan-Margulis relative property (T). Let $A$ be a discrete abelian group and let $\Gamma$ be a subgroup of $\textrm{Aut}(A)$. Let $\mu$ be the Haar measure on $\hat{A}$, the dual group of $A$. Then, the action $\Gamma \curvearrowright (\hat{A},\mu)$ has the property (T) relative to the space if and only if the pair of groups $(A\rtimes\Gamma, A)$ has the relative property (T) (see Proposition 5.1 \cite{Pop}). 

\medskip

Consider in particular the following pair of groups $(\Z^2\rtimes\Gamma,\Z^2)$ with the natural action of $\Gamma\subset$SL$_2(\Z)$ on $\Z^2$. Burger \cite{Bur} showed that if $\Gamma$ is non-amenable, then the pair $(\Z^2\rtimes\Gamma,\Z^2)$ has the relative property (T). In other words, the action $\Gamma\curvearrowright (\mathbb{T}^2,\lambda)$ has the property (T) relative to the space when $\Gamma$ is non-amenable.

\medskip

Property (T) relative to the space have led to several remarkable applications; in Popa's work on fundamental groups of von Neumann $II_1$ factors \cite{Pop}, in the constructions of orbit inequivalent actions of non-amenable groups \cite{GabPop} \cite{Ioa2} \cite{Eps}, to produce actions of $\mathbb{F}_\infty$ (the free group on infinitely many generators) whose $II_1$ factors and orbit equivalence relations have prescribed fundamental group (\cite{PopVa}).

\medskip

Yet, few examples of groups admitting actions with the property (T) relative to the space are known. Since probability measure preserving actions of amenable groups never have the property (T) relative to the space, Popa raised the following natural question (see Problem 5.10.2 \cite{Pop}): does every non-amenable countable group admit a free (i.e. essentially free) ergodic action with the property (T) relative to the space?

\medskip
This question is still open. However some partial results were obtained.

\medskip

It follows from Theorem 1.1 in \cite{Fer} that every finitely generated group $\Gamma$ for which there exists a morphism $\varphi:\Gamma\rightarrow$ SL$_n(\R)$ such that the $\R$-Zariski-closure $\overline{\varphi(\Gamma)}^Z(\R)$ is non-amenable has an action $\Gamma\curvearrowright (X,\mu)$ with the property (T) relative to the space.

Besides, it has been showed in \cite{Gab}, Theorem 1.3, that every non-amenable free product of groups admits such actions. 

In addition, it was proved in \cite{Ioa2}, Theorem 4.3, that every non-amenable group admits free ergodic measure preserving actions with a weak form of the property (T) relative to the space.

\medskip

Most of the actions with this property were built using the relative property (T) of a pair $(A\rtimes\Gamma, A)$, where $A$ is a discrete abelian group and $\Gamma$ a subgroup of $\textrm{Aut}(A)$.

\medskip

In this paper, we study the case of finitely generated linear groups and we show that we can get rid of the assumption made by Fernos.

\begin{thm}
Every finitely generated non-amenable linear group $\Gamma$ over a field of characteristic zero admits an ergodic action $\Gamma\curvearrowright (X,\mu)$ with the property (T) relative to the space.
\label{thml}
\end{thm}

The actions we construct usually have a kernel and therefore are not necessarily free. 
Still, in the case where $\Gamma$ has trivial solvable radical, we are able to produce free rigid actions.

\begin{thm}
Every finitely generated non-amenable linear group $\Gamma$ over a field of characteristic zero which has trival solvable radical admits a free ergodic action $\Gamma\curvearrowright (X,\mu)$ with the property (T) relative to the space.

\label{thmr}
\end{thm}

\medskip

In order to obtain these theorems, we had to investigate natural questions concerning the stability of the property (T) relative to the space through standard constructions as follows.

\medskip

Let $\Gamma$ be a countable group acting on $(X,\mu)$ by measure preserving transformations. Let $\Lambda$ be a subgroup of $\Gamma$. If the restricted action $\Lambda\curvearrowright (X,\mu)$ has the property (T) relative to the space, then obviously the action $\Gamma\curvearrowright (X,\mu)$ does as well.

\medskip

We show that the converse holds when $\Lambda$ is a co-amenable subgroup of $\Gamma$ (see Definition \ref{defEym}).

\begin{prop}
Let $\Gamma$ be a countable group and $\Lambda$ a co-amenable subgroup.

A probability measure preserving action $\Gamma\curvearrowright (X,\mu)$ has the property (T) relative to the space if and only if the restricted action $\Lambda\curvearrowright (X,\mu)$ has the property (T) relative to the space.

\label{lemM}
\end{prop}

\medskip

Another natural question concerns the stability of rigid actions by co-induction (see section 2.5 for the definition). Let $\Gamma_0\subset \Gamma$ be countable groups. Let $\Gamma_0\curvearrowright (X,\mu)$ be a probability measure preserving action with the property (T) relative to the space.

\begin{prop}
The co-induced action $\Gamma\curvearrowright (Y,\eta)$ has the property (T) relative to the space if and only if $[\Gamma:\Gamma_0]<\infty$.

\label{lemcoind}
\end{prop}

\medskip

The crucial tool for the proof of our results is an extension of our characterization of rigid actions on homogeneous spaces of real Lie groups (Theorem 1.3, \cite{Moi}) to the framework of $p$-adic groups, where $p\in\mathcal{P}$, the set of prime numbers.

%
%

\medskip

We need the following notations. Let $G$ be a $p$-adic or real Lie group and let $\Lambda$ be a lattice in $G$. Let Aut$(G)$ be the group of continuous automorphisms of $G$ and Aut$_\Lambda(G)$ the subgroup of Aut$(G)$ of all $\sigma\in $Aut$(G)$ with $\sigma(\Lambda)=\Lambda$. Set Aff$(G)=G\rtimes$Aut$(G)$ and $$\mathrm{Aff}_\Lambda(G)=G\rtimes\mathrm{Aut}_\Lambda(G).$$

The group of affine transformations of $G/\Lambda$, $\mathrm{Aff}_\Lambda(G)$, acts in a natural way on $G/\Lambda$. Any $G$-invariant probability measure on $G/\Lambda$ is also invariant under Aff$_\Lambda(G)$.

Moreover, the linear action of Aff$(G)$ on $\frak{g}$ induces an action on the projective space $\mathbb{P}(\frak{g})$.

\begin{thm}
Let $G$ be a $p$-adic or real Lie group. Let $\Lambda\subset G$ be a lattice and let $\mu$ be the $G$-invariant probability measure on the corresponding homogeneous space $G/\Lambda$.
Let $\Gamma$ be a countable subgroup of $\mathrm{Aff}_\Lambda(G)$. 
Then the following properties are equivalent :

\begin{itemize}
\item[(i)] $\Gamma \curvearrowright (G/\Lambda,\mu)$ has the property (T) relative to the space;

\item[(ii)] there is no $\Gamma$-invariant probability measure on $\mathbb{P}(\frak{g})$.

\end{itemize} 
\label{thmp}

\end{thm}

Note that by Lemma \ref{lemPr}, this Theorem extends to the framework of products of $p$-adic and real Lie groups.

\medskip

Let $S$ be a finite subset $\mathcal{P}\cup\{\infty\}$. Let $G$ be a group of the form $$\prod_{p\in S}G_p$$where for every $p\in S$, $G_p$ is a $p$-adic Lie group, and $G_\infty$ is a real Lie group, if $\infty\in S$. The Lie algebra of $G$ is$$\frak{g} := \oplus_{p\in S} \frak{g}_p$$the direct sum of the Lie algebras $\frak{g}_p$ of $G_p$. Let $\Lambda\subset G$ be a subset of the form $$\prod_{p\in S}\Lambda$$where for every $p\in S$, $\Lambda_p$ is a lattice of $G_p$. Let $\mu$ be the $G$-invariant probability measure on the homogeneous space $G/\Lambda$.

\medskip 

In this setting, by Theorem \ref{thmp} and Lemma \ref{lemPr} we deduce the following result.

\begin{thm}
Let $\Gamma$ be a countable subgroup of $\prod_{p\in S}\mathrm{Aff}_{\Lambda_p}(G_p)$. 

The following assertions are equivalent:

\begin{itemize}
\item[(i)] $\Gamma \curvearrowright (G/\Lambda,\mu)$ has the property (T) relative to the space;

\item[(ii)] for every $p\in S$, there is $\Gamma$-invariant probability measure on $\mathbb{P}(\frak{g}_p)$.

\end{itemize}
\label{TheS}
\end{thm}

\medskip

We now explain how the actions arise in our Theorems \ref{thml} and \ref{thmr}; we construct actions by translations on homogeneous spaces $(G/\Lambda,\mu)$ so that there is no $\Gamma$-invariant probability measure on $\mathbb{P}(\frak{g})$. In this case, Theorem \ref{thmp} applies.

\medskip

If $\Gamma\subset G$ is an unbounded Zariski-dense subgroup of a simple Lie group $G$, then there is no $\Gamma$-invariant probability measure on $\mathbb{P}(\frak{g})$. Indeed, by Furstenberg's Lemma, we would get a finite index subroup $\Gamma_0\subset\Gamma$ which would stabilise a proper subspace $V\subset \frak{g}$. As such, we would obtain an proper ideal in $\frak{g}$, which would contradict the simplicity of $G$.

\medskip

We show that every non-amenable finitely generated linear group admits such a representation.

\begin{prop}
Let $\Gamma\subset$ GL$_n(k)$ be a finitely generated group over a field of charasteristic zero. Assume that Rad$(\Gamma)=\{e\}$. Then there exist local fields $k_1,\dots,k_r$ of characteristic zero, and simple Lie groups $S_1,\dots,S_r$ over $k_1,\dots,k_r$ with an injective morphism $$j:\Gamma\hookrightarrow S_{1}\times\dots\times S_{r}$$ such that the projection of $j(\Gamma)$ on each factor $S_i$ is unbounded and Zakiski-dense.
\label{lemdecom}
\end{prop}

\medskip

Rigidity for actions $\Gamma\curvearrowright (G/\Lambda,\mu)$ by translations was first investigated by Ioana and Shalom \cite{Ioa} in the case where $G$ is real algebraic group. They observed that $\Gamma\curvearrowright (G/\Lambda,\mu)$ is not rigid if $\Gamma$ has infinite centre. Consequently they asked whether $\mathbb{F}_2\times\Z$ admits a free ergodic rigid action at all. 

\medskip

We show (see Proposition \ref{propmoi}) that indeed $\Gamma\times\Z$ has a free ergodic action with the property (T) relative to the space for every finitely generated linear group $\Gamma$ with trivial solvable radical satisfying Fernos' condition that is: $\Gamma$ admist a morphism $\varphi:\Gamma\rightarrow$ SL$_n(\R)$ such that the $\R$-Zariski-closure $\overline{\varphi(\Gamma)}^Z(\R)$ is non-amenable.

\medskip

Note that by Proposition \ref{lemM}, it is enough to find a free ergodic action $\Gamma\curvearrowright (X,\mu)$ which has the property (T) relative to the space and such that the commutant of $\Gamma$ in Aut$(X,\mu)$  has an element of infinite order (where Aut$(X,\mu)$ is the set of Borelian bijections which preserve the measure $\mu$, identified up to null set).

\medskip

For such a group, one can build an action by automorphisms on the dual group $(A,\lambda)$ of a certain $\Z[S^{-1}]^n$ which has the property (T) relative to the space, where $S$ is some finite set of rational primes. The diagonal action of $\Gamma$ on $(A\times A,\lambda\otimes \lambda)$ still has the property (T) relative to the space. Then, we easily obtain a free action of $\Z$ on $(A\times A,\lambda\otimes \lambda)$ which commutes with $\Gamma$. 

\medskip

As mentioned before, the action of $\Gamma$ on $(A\times A,\lambda\otimes \lambda)$ usually has a kernel. Therefore we apply Theorem \ref{thmr} to get a free ergodic action of $\Gamma$ which has the property (T) relative to the space. Considering these two actions, we finally obtain a free ergodic action of $\Gamma\times \Z$ which has the property (T) relative to the space.

\medskip

The paper is organized as follows. Section 2 contains a few preliminary results, among which we show that the property (T) relative to the space is stable for co-amenable subgroups. The proof of Theorem \ref{thmp} is given in Section 3. In Section 4, we prove the theorems \ref{thml} and \ref{thmr}, and we detail the answer to Ioana-Shalom's question.

\bigskip

\textbf{Acknowledgments.} 
I am very grateful to Bachir Bekka for valuable discussions concerning this work, and for his helpful advices. I also owe many thanks to Damien Gaboriau for his useful comments.

\section{Preliminaries}




\subsection{Ergodic characterization of the property (T) relative to the space}

Ioana gave the following purely ergodic theoretical characterization of the property (T) relative to the space (see Theorem 4.4 \cite{Ioa3}, Proposition 1 \cite{Ioa}). We will use this criteria instead of the original characterization given in terms of pair of von Neumann algebras. We denote by $B(X)$ the algebra of bounded measurable complex valued functions on $X$.

\begin{prop}[\cite{Ioa}]
\label{propIoa}
A measure preserving action $\Gamma\curvearrowright (X,\mu)$ of a countable group $\Gamma$ on a probability space $(X,\mu)$ has the property (T) relative to the space if and only if, for any sequence of Borel probability measures $\nu_n$ on $X\times X$ satisfying:
\begin{enumerate}
\item $p_*^i\nu_n=\mu$ for all $n$ and $i=1,2$, where $p^i:X\times X\rightarrow X$ denotes the projection onto the i-th coordinate,
\item $\int_{X\times X} \varphi(x)\psi(y)d\nu_n(x,y)\underset{n\rightarrow\infty}{\longrightarrow}\int_X \varphi(x)\psi(x)d\mu(x)$, for all bounded Borel functions $\varphi,\psi\in B(X)$,
\item $||(\gamma\times \gamma)_*\nu_n-\nu_n||\underset{n\rightarrow\infty}{\longrightarrow}0$, for every $\gamma\in\Gamma,$
\end{enumerate} we have that $\nu_n(\Delta_X) \underset{n\rightarrow\infty}{\longrightarrow} 1$ where $\Delta_X$ denotes the diagonal in $X\times X$.
\end{prop}

\subsection{Diagonal action on a finite product of standard probability spaces}

We now study the behavior of the property (T) relative to the space by passing to diagonal actions on product spaces.

\begin{lem}
For $n\in \N$, let $(X_i,\mu_i)_{i=1,\dots,n}$ be standard probability spaces. Set $(X,\mu)=\prod_i(X_i,\mu_i)$, and Aut$(X,\mu)$ is the set of Borelian bijections which preserve the measure $\mu$, identified up to null set. Let $\Gamma\subset \prod_i $ Aut$(X_i,\mu_i)$ be a countable group.

Then the following properties are equivalent:

\begin{itemize}
\item[(i)] $\Gamma\curvearrowright (X,\mu)$ has the property (T) relative to the space,
\item[(ii)] for every $i=1,\dots,n,$ the action $\pi_i(\Gamma)\curvearrowright (X_i,\mu_i)$ has the property (T) relative to the space, where $\pi_i$ is the projection of $\prod_{j=1}^n $ Aut$(X_j,\mu_j)$ onto Aut$(X_i,\mu_i)$. 
\end{itemize}
\label{lemPr}
\end{lem} 

\textbf{Proof:}

\textbf{\underline{(i) $\Rightarrow$ (ii)}:}

Let $i\in\{1,\dots,n\}$. Let $\nu_n$ be a sequence of probability measures on $X_i\times X_i$ with the following properties: 

\begin{enumerate}
\item $p_*^i\nu_n=\mu_i$, for all $n$ and $i=1,2,$
\item $\int_{X_i\times X_i}f(x)g(y)\mathrm{d}\nu_n(x,y)\underset{n\rightarrow\infty}{\longrightarrow} \int_{X_i} f(x)g(x)\mathrm{d}\mu_i$, for all $f,g\in\mathrm{B}(X_i),$
\item $\parallel (\pi_i(\gamma)\times\pi_i({\gamma})_* \nu_n -\nu_n \parallel \underset{n\rightarrow\infty}{\longrightarrow} 0$, for all $\gamma\in \Gamma.$
\end{enumerate}
We have to show that $\underset{n\rightarrow \infty}{\lim}\, \nu_n(\Delta_{X_i}) = 1.$ Let $\eta_n\in \mathcal{M}(X_i\times X_i)$ be the sequence of probability measures defined by: 
\[\eta_n=\nu_n\otimes(\bigotimes_{j\neq i}q_*^j\mu_j),\]where $q^j:X_j\rightarrow X_j\times X_j,\,x\mapsto(x,x)$. This sequence satisfies the conditions 1. 2. and 3. of Proposition \ref{propIoa}. 

Since the action $\Gamma\curvearrowright (X,\mu)$ has the property (T) relative to the space, we have:
\[\eta_n(\Delta_X)=\nu_n(\Delta_{X_i})\underset{n\rightarrow\infty}{\longrightarrow}1.\] Therefore, the action $\pi_i(\Gamma)\curvearrowright (X_i,\mu_i)$ has the property (T) relative to the space.

\bigskip

\textbf{\underline{(ii) $\Rightarrow$ (ii)}:}

Assume now that for $i=1,\dots,n$, the actions $\pi_i(\Gamma)\curvearrowright (X_i,\mu_i)$ have the property (T) relative to the space. Let $\eta_n$ be a sequence of probability measures on $X\times X$ satisfying the three conditions of Proposition \ref{propIoa}. We have to show that $\underset{n\rightarrow \infty}{\lim}\, \eta_n(\Delta_{X}) = 1.$

For $i=1,\dots,n$, let $\nu_n^i$ be the sequence of probability measures on $(X_i\times X_i,\mu_i\otimes \mu_i)$ defined by:
\[\nu_n^i=p_*^i\eta_n,\,\text{where }p_*^i \text{ is the projection } p^i:X\times X\rightarrow X_i\times X_i .\]

For all $i=1,\dots,n$, the sequence $\nu_n^i$ satisfies the properties of Proposition \ref{propIoa} for the action $\pi_i(\Gamma)\curvearrowright (X_i,\mu_i)$. Hence, for all $i=1,\dots,n$:
\[\nu_n^i(\Delta_{X_i})=\eta_n( \prod_{j\neq i} (X_j\times X_j)\times \Delta_{X_i})\underset{n\rightarrow\infty}{\longrightarrow} 1.\]Therefore,
\[ \eta_n(\Delta_X)\underset{n\rightarrow\infty}{\longrightarrow} 1,\] and the action $\Gamma\curvearrowright (X,\mu)$ has the property (T) relative to the space.

\begin{flushright}
$\square$
\end{flushright}

\begin{rk}

The result in the previous lemma does not hold for an infinite product $(\prod_{i=1}^\infty X_i,\mu^{\otimes \N})$ of standard probability spaces $(X_i,\mu_i)$. Indeed, the following sequence of measures:
\[\eta_n =(\bigotimes_{i=1}^n q_* \mu ) \otimes (\bigotimes_{i=n+1}^\infty  \mu\otimes \mu)\] where $q: X \rightarrow X \times X, x \mapsto (x,x)$ satisfies properties 1. 2. and 3. of Proposition \ref{propIoa}. Since $\eta_n(\Delta_X)=0$ for all $n$, the diagonal action on a infinite product of standard probability spaces never has the property (T) relative to the space
\end{rk}

\subsection{Proof of Proposition \ref{lemM}}

Recall the notion of a co-amenable subgroup or an amenable homogeneous space in Eymard's sense \cite{Eym}.
\begin{defi}
Let $\Gamma$ be a countable group and let $\Lambda$ be a subgroup of $\Gamma$. We say that $\Lambda$ is a co-amenable subgroup of $\Gamma$ if there exists a sequence $f_n\in l^1(\Gamma/\Lambda)$ such that 
\begin{enumerate}
\item $f_n \geq 0$
\item $||f_n||_1=1$
\item $||\pi(\gamma)f_n-f_n||_1 \underset{n\rightarrow\infty}{\longrightarrow} 0$ for all $\gamma\in \Gamma$, where $\pi(\gamma)f_n(x)=f_n(\gamma^{-1}x)$, for all $x\in \Gamma/\Lambda$.
\end{enumerate}
\label{defEym}
\end{defi}

We now proceed with the proof of Proposition \ref{lemM}.

Let $D\subset \Gamma$ be a set of representatives for the left cosets of $\Lambda$ in $\Gamma$ so that $$\Gamma=\sqcup_{\lambda\in\Lambda} D\lambda.$$ We transfer an action of $\Gamma$ on $\Gamma/\Lambda$ to an action of $\Gamma$ on $D$: for all $\gamma\in \Gamma$, $h\in D$, set

\[\gamma\tilde . h := \gamma h c(\gamma,h)^{-1}, \]$\text{where } c (\, ,\, ) : \Gamma\times D \rightarrow \Lambda\text{ is the cocycle defined by }  \gamma h c(\gamma,h)^{-1}\in D.$

Since $\Gamma/\Lambda$ is amenable, there exists $\tilde{f_n}\in l^1(\Gamma/\Lambda)$ a sequence of non-negative functions, such that $||f_n||_1=1$ and for all $\gamma\in \Gamma$, $||\pi(\gamma)\tilde{f_n}-\tilde{f_n}||_1 \underset{n\rightarrow\infty}{\longrightarrow} 0 $. Let $f_n\in l^1(D)$ be the sequence of functions such that

\[ f_n(h)=\tilde{f_n}(h\Lambda),\quad\text{for all }h\in D.\]




\medskip

Assume, by contrapositive, that the action of $\Lambda$ on $(X,\mu)$ is not rigid. 

Let $\nu_n$ be a sequence of Borel probability measures on $X\times X$ with the following properties: 

\begin{enumerate}
\item $p_*^i\nu_n=\mu$ for all $n$ and $i=1,2$,
\item $\int_{X\times X} \phi(x)\psi(y)d\nu_n(x,y) \underset{n\rightarrow\infty}{\longrightarrow}\int_X \phi(x)\psi(x)d\mu(x)$, for all bounded Borelian functions $\phi,\psi\in B(X)$,
\item $||(\lambda \times \lambda)_*\nu_n-\nu_n||\underset{n\rightarrow\infty}{\longrightarrow}0$, for all $\lambda\in \Lambda$,
\item $\underset{n\rightarrow\infty}{\varliminf} \nu_n(\Delta_X)< 1.$ 
\end{enumerate}

Let us consider the set of measures $\eta_{n,m}$ defined on $X\times X$ by :

\[\eta_{n,m} = \sum_{h\in D} f_m(h) (h\times h)_*\nu_n \]

We claim that there exists a sequence $\eta_m\in\{\eta_{n,m}\,:\, (n,m)\in \N\times \N\}$ which satisfies the previous four properties for the action of $\Gamma$ on $(X,\mu)$. 

\medskip
As for the first point, one directly checks that for all $n,m \in \N$ and $i=1,2$, we have $p^i_*\eta_{n,m}=\mu$.

\medskip
We now look at the second property.

Let $\phi,\psi\in $ B$(X)$ be two bounded Borel functions. 

\[
\int_{ X\times X} \phi(x)\psi(y)\mathrm{d}\eta_{n,m}(x,y)=\sum_{h\in D}f_m(h) \int_{ X\times X} \phi(hx)\psi(hy)\mathrm{d}\nu_n(x,y).
\]

But for all $h\in D\subset \Gamma,$ one has
\[  \int_{ X\times X} \phi(hx)\psi(hy)\mathrm{d}\nu_n(x,y)  \underset{n\rightarrow\infty}{\longrightarrow} \int_{ X} \phi(x)\psi(x)\mathrm{d}\mu(x),\]
by the $\Gamma$-invariance of the measure $\mu$.

Hence, we get for all $m\in \N$

\[ \int_{ X\times X} \phi(x)\psi(y)\mathrm{d}\eta_{n,m}(x,y) \underset{n\rightarrow\infty}{\longrightarrow} \int_{ X} \phi(x)\psi(x)\mathrm{d}\mu(x).\]

Therefore, the second property as soon as $n$ tends to infinity for fixed $m$.

\medskip

We will now focus on the third point. Let $\gamma\in\Gamma$. For all $\phi \in $ B$(X\times X)$ with $||\phi||_\infty=1$, we have:
 
\[|\int_{x\in X\times X} \phi(\gamma.x)-\phi(x) \mathrm{d}\eta_{n,m}(x)| = |\int_{x\in X\times X} \sum_{h\in D}  (\phi (\gamma.(h.x))-\phi(h.x) )f_m(h) \mathrm{d}\nu_n(x)| .\]

But, one can write:

\[ \sum_{h\in D} \phi(\gamma.( h.x))f_m(h) = \sum_{h\in D} \phi((\gamma\tilde{.}h).( c(\gamma,h) .x))f_m(h)\]

with $c(\gamma,h)\in \Lambda$. Then by substitution, 

\[ \sum_{h\in D} \phi(\gamma.( h.x))f_m(h) = \sum_{h\in D} \phi(h c(\gamma,\gamma^{-1}\tilde{.}h) .x)f_m(\gamma^{-1}\tilde{.}h).\]

Finally, we get

\[\begin{split}
&|\int_{x\in X\times X} \phi(\gamma.x)-\phi(x) \mathrm{d}\eta_{n,m}(x)|\\ =& |\int_{x\in X\times X} \sum_{h\in D} \phi(h c(\gamma,\gamma\tilde{.}h) .x)f_m(\gamma^{-1}\tilde{.}h)-\phi(h.x)f_m(h)  \mathrm{d}\nu_n(x)| \end{split}.\]

It can also be written 

\[\begin{split}
&|\int_{x\in X\times X} \phi(\gamma.x)-\phi(x) \mathrm{d}\eta_{n,m}(x)| \\
\leq& |\int_{x\in X\times X} \sum_{h\in D} ( \phi(h c(\gamma,\gamma^{-1}\tilde{.}h) .x)f_m(\gamma^{-1}\tilde{.}h)-\phi(hc(\gamma,\gamma^{-1}\tilde{.}h) .x)f_m(h))  \mathrm{d}\nu_n(x)| \\
&+ |\int_{x\in X\times X} \sum_{h\in D}( \phi(h c(\gamma,\gamma^{-1}\tilde{.}h) .x)f_m(h) -\phi(h.x) f_m(h)) \mathrm{d}\nu_n(x)|, \end{split}\]

where the first integral on the right hand side is

\[\begin{split}
&|\int_{x\in X\times X} \sum_{h\in D} ( \phi(h c(\gamma,\gamma^{-1}\tilde{.}h) .x)f_m(\gamma^{-1}\tilde{.}h)-\phi(hc(\gamma,\gamma^{-1}\tilde{.}h) .x)f_m(h))  \mathrm{d}\nu_n(x)|\\
\leq&   ||\pi(\gamma)f_n-f_n||_1,
\end{split}\]

and the second one is

\[\begin{split}|\int_{x\in X\times X} \sum_{h\in D}&( \phi(h c(\gamma,\gamma^{-1}\tilde{.}h) .x)f_m(h) -\phi(h.x) f_m(h))  \mathrm{d}\nu_n(x)| \\
=& |\int_{x\in X\times X} \sum_{h\in D}( \phi(h c(\gamma,\gamma^{-1}\tilde{.}h) .x) -\phi(h.x) )f_m(h)  \mathrm{d}\nu_n(x)|\\
\leq& \sum_{h\in D} ||(k_{\gamma,h}\times k_{\gamma,h})_*\nu_n -\nu_n || f_m(h),
\end{split}
\]
where $k_{\gamma,h}=c(\gamma,\gamma^{-1}\tilde{.}h)\in \Lambda$. Hence, we get: 

\[ || (\gamma\times \gamma)_* \eta_{n,m} - \eta_{n,m} || \leq ||\pi(\gamma)f_m-f_m||_1 + \sum_{h\in D} ||(k_{\gamma,h}\times k_{\gamma,h})_*\nu_n -\nu_n || f_m(h).\]

The first term on the right hand side tends to zero as $m$ tends to infinity. Moreover, the second term tends to zero, for $m$ fixed. Hence, we can easily build a sequence satisfying the second and the third property. Let us call this subsequence $\eta_m$.

\medskip

Finally, since for all $n,m\in N$ we have $\eta_{n,m}(\Delta_X)=\nu_n(\Delta_X)$, the fourth property is also verified for our sequence $\eta_m$.
\medskip

Therefore, we get a sequence $\eta_m$ of measures on $X\times X$ showing that the action of $\Gamma$ on $(X,\mu)$ does not have the property (T) relative to the space.

\begin{flushright}
$\square$
\end{flushright}

\subsection{Co-induction}

We first recall the construction of the co-induction action.

Let $\Gamma_0\subset \Gamma$ be countable groups. Let $\Gamma_0\curvearrowright (X,\mu)$ be a probability measure preserving action.

Let $D$ be a fundamental domain for the right action of $\Gamma_0$ on $\Gamma$. Let $(Y,\eta)=(\prod_D X,\otimes_D \mu)$ and consider the action of $\Gamma$ on $(Y,\eta)$ given by:

\[\gamma y=\gamma (x_d)_d=(\sigma(\gamma^{-1},d)^{-1} x_{\gamma^{-1}.d})_{d}\] where the action of $\Gamma$ on $D$; $\Gamma\times D \rightarrow D, (\gamma,d)\mapsto \gamma . d$ is given by the cocycle $\sigma:\Gamma\times D\rightarrow \Gamma_0$.

\textbf{Proof of Proposition \ref{lemcoind}:}

For the direct implication, we proceed by contradiction. Assume that $[\Gamma:\Gamma_0]=\infty$.

Consider the sequence of measures $\mu_n$ defined on $X\times X$ by: $$\nu_n:= (1-\frac1{n})q_*\mu + \frac1{n} \mu\otimes \mu,$$ where $q: X \mapsto X\times X, x\mapsto (x,x) $. 

This sequence of measures is invariant by the diagonal action of $\Gamma_0$ and satisies $\nu_n(\Delta_X)=1-\frac1{n}$.

We now introduce the sequence of measures $eta_n$ defined on $Y\times Y$ by: $$\eta_n =\bigotimes_{D} \nu_n.$$

This sequence of measures satisfies the properties of \ref{propIoa}. Yet, we have for all $n\in\N$, $$\eta_n(\Delta_Y)=0,$$ since the product is infinite.

This contradicts the fact the co-indued action has the property (T) relative to the space.

Therefore, $[\Gamma:\Gamma_0]<\infty$.

\medskip

Assume now that the action $\Gamma_0\curvearrowright (X,\mu)$ has the property (T) relative to the space and $[\Gamma:\Gamma_0]<\infty$.

Set $m=[\Gamma:\Gamma_0]\in\N$ and $D=\{ d_1,\dots,d_m \}$.

Let $\eta_n\in\mathcal{M}(Y\times Y)$ be a sequence of probability measures satisfying the properties of Proposition \ref{propIoa}. It is enough to show that for $i=1,\dots,q$ we have $q^i_*\eta_n(\Delta_X)\rightarrow 1$, where $q^i:Y\times Y \rightarrow X\times X$ is the projection onto the $d_i$ coordinate.

\medskip

The sequence of measures $q^i_*\eta_n\in\mathcal{M}(X \times X)$ satisfies properties 1) and 2) of Propostion \ref{propIoa}. Since $\Gamma_0\curvearrowright (X,\mu)$ has the property (T) relative to the space, it sufficies to show that this sequence satisfies also property 3). 

\medskip

Let $\gamma_0\in\Gamma_0$. Let $\varphi\in $ B$(X\times X)$. For any $n\in N$, we have:
\[ 
\begin{split}
(\gamma_0\times\gamma_0)_*(q^i_*\eta_n)(\varphi)-q^i_*\eta_n(\varphi)&=\int_{x\in X\times X}( \varphi(\gamma_0.x)-\varphi(x))\mathrm{d}q^i_*\eta_n(x)\\
&=\int_{(x_1,\dots,x,\dots,x_m)\in Y\times Y}(\varphi(\gamma_0.x)-\varphi(x))\mathrm{d}\eta_n(x_1,\dots,x,\dots,x_q)\\
&=\int_{y\in Y\times Y} \psi(\gamma.y)-\psi(y)\mathrm{d}\eta_n (y)
\end{split}
\]where $x$ is on the i-th coordinate, $\psi(y)=\psi(x_1,\dots,x,\dots,x_m)=\varphi(x)$ and $\gamma=d_i \gamma_0^{-1} d_i^{-1}$. As such, we get:
\[
\begin{split}
|(\gamma_0\times\gamma_0)_*(q^i_*\eta_n)(\varphi)-q^i_*\eta_n(\varphi)|&= |(\gamma\times \gamma)_*\eta_n(\psi)-\eta_n(\psi)|\\
&\leq  ||(\gamma\times \gamma)_*\eta_n-\eta_n||.
\end{split}
\]Hence $q^i_*\eta_n$ satisfies the properties of Proposition \ref{propIoa}. Therefore $q^i_*\eta_n(\Delta_X)\rightarrow 1$ for $i=1,\dots,m$ and the co-induced action has the property (T) relative to the space.

\begin{flushright}
$\square$
\end{flushright}


\section{Rigidity for group actions on homogeneous spaces}

\subsection{Notations}

Let $G$ be a $p$-adic or real Lie group, let $\frak{g}$ be the Lie algebra of $G$ and denote by $p : \frak{g}  \backslash \{ 0\} \rightarrow \mathbb{P}(\frak{g})$ the canonical projection onto the projective space $\mathbb{P}(\frak{g})$ of $\frak{g}$. 

Let $q: G \rightarrow \frak{g}$ be any Borel map which equals the logarithm in some neighborhood of $e \in G$. Define $r : (G\times G)  \backslash  \Delta \rightarrow G,\,(x,y)\mapsto xy^{-1}$. 

As in $\cite{Ioa}$,we consider the map $$\rho\, : \,(G\times G)  \backslash  \Delta \rightarrow \mathbb{P}(\frak{g}),\,\rho(x,y)=p(q(r(x,y))).$$ We will denote by Ad$: G \rightarrow GL(\frak{g})$ the adjoint representation of $G$ on $\frak{g}$ as well as the associated action $G\rightarrow $Aut$(\mathbb{P}(\frak{g}))$ of $G$ on $\mathbb{P}(\frak{g})$.

Recall that the affine group Aff$(G)$ is the semi-direct product $G\rtimes$ Aut$(G)$. The associated action of Aff$(G)$ on $\mathbb{P}(\frak{g})$ is given for all $\gamma=(g,\sigma)\in\text{Aff}(G)$ and for all $Y\in\mathbb{P}(\frak{g})$ by: $$\gamma.\bar{Y}=\mathrm{Ad}(g)\circ\mathrm{d}\sigma_e(\bar{Y}),$$ where $\mathrm{d}\sigma_e$ is the derivated automorphism at the identity of $\sigma$.

\subsection{Proof of Theorem \ref{thmp}}

We proceed by contrapositive by showing that the following properties are equivalent:
\begin{itemize}
\item[(i')] $\Gamma \curvearrowright (G/\Lambda,\mu)$ does not have the property (T) relative to the space,
\item[(ii')] there exists a $\Gamma$-invariant probability measure on $\mathbb{P}(\frak{g})$.

\end{itemize} 

\medskip

\textbf{\underline{(i') $\Rightarrow$ (ii')}:}  the proof is similar to the proof of Theorem 1.3 in \cite{Moi} and will be omitted.

\bigskip

\textbf{\underline{(ii') $\Rightarrow$ (i')}:} Assume now that there exists a $\Gamma$-invariant probability measure $\zeta$ on $\mathbb{P}(\frak{g})$.

\medskip

By Proposition \ref{lemM}, we know that the action of $\Gamma$ on $(X,\mu)$ has the property (T) relative to the space if, and only if, the action of $[\Gamma,\Gamma]$ on $(X,\mu)$ has the property (T) relative to the space. Hence, we can assume that Ad$(\Gamma)\subset$ SL$(\frak{g})$. 

We have to consider two cases as to whether the closure of Ad$(\Gamma)$  in the Hausdorff topology of SL$(\frak{g})$ is compact or not.

\medskip

\emph{First case:} Assume that $K:=\overline{\Gamma}$ is compact.

Let $\lambda$ be the normalized Haar measure on $K$ and let $Y\in\frak{g}\backslash\{0\}$ be in a neighborhood of $0$ such that for all $n\in\N^*$, $x_n=\exp(Y/n)$ is well defined.

Consider the sequence of probability measures $\eta_n$ on $X\times X$ defined by:
\[ \eta_n(\varphi)=\int_K\int_X\varphi(x,k(x_n)x)\mathrm{d}\mu(x)\mathrm{d}\lambda(k),\quad\text{for all } \varphi\in\text{ B}(X\times X),\]where $\lim_{n\rightarrow\infty} k(x_n):=\exp(k(Y/n))$. Observe that, $\lim_{n\rightarrow\infty} =e$ and that $k(x_n)\neq e$ for all $n\geq 1$. 

The sequence $\eta_n$ satisfies the three properties of Proposition \ref{propIoa} ($\eta_n$ is even $\Gamma$-invariant). Yet, $\eta_n(\Delta_X)=0$ for all $n$. Hence, the action does not have the property (T) relative to the space.

\medskip

\emph{Second case:} Assume that $\Gamma\subset$ SL$(\frak{g})$ is unbounded.

Since $\zeta$ is a $\Gamma$-invariant measure on $\mathbb{P}(\frak{g})$, by Furstenberg's Lemma, there exists a finite index subgroup $\Gamma_0$ of $\Gamma$ which fixes a proper non-zero subspace $V_0\subset \frak{g}$.

\medskip
Again, we have to distinguish two cases. If $\Gamma_0\subset $ GL$(V_0)$ is relatively compact, then we can apply our first case to $K_0:=\overline{\Gamma_0}\subset$ GL$(V_0)$. So, the action of $\Gamma_0$ on $X$ does not have the property (T) relative to the space. Proposition \ref{lemM} implies that the action of $\Gamma$ on $(X,\mu)$ does not have the property (T) relative to the space.

\medskip

Assume now that $\Gamma_0\subset $ SL$(V_0)$ is unbounded. By Furstenberg's Lemma, there exists a proper non-zero subspace $W_1\subset V_0$ which is invariant under a finite index subgroup $\Gamma_1$ of $\Gamma$, and such that the invariant measure $\zeta$ is supported by the image of $W_1$ in $\mathbb{P}(V_0)$. We may now consider the restriction of the $\Gamma_1$-action to $W_1$.

\medskip

Continuing this way, we obtain a proper subspace $W$ which is invariant under a subgroup of finite index of $\Gamma$ such that either the image of $\Gamma$ in GL$(W)$ is bounded or $W$ is a line.

\medskip

If the image is bounded, we conclude as in the first case that $\Gamma\curvearrowright (X,\mu)$ does not have the property (T) relative to the space. So, we may assume that $W$ is a line.

Then, $[\Gamma_1,\Gamma_1]$ has an non-trivial fixed point $Y$ in $\frak{g}$, and the action $[\Gamma_1,\Gamma_1]\curvearrowright (X,\mu)$ does not have the property (T) relative to the space.

Indeed, we can consider, for $n$ large enough, the sequence given by $x_n=\exp(Y/n)$. Assume that this sequence is well defined for $n\geq 1$. Then $x_n\in G$ is fixed by $[\Gamma_1,\Gamma_1]$. Since $x_n\rightarrow e$ and $x_n\neq e$ for $n$ large enough, and since $\Lambda$ is discrete, we can assume that $x_n\notin \Lambda$ for every $n\geq 1$.

Consider the map
\[
\begin{split}
p_n:G/\Lambda & \rightarrow G/\Lambda \times G/\Lambda \\
x & \mapsto (x , x_n x)
\end{split}
\] for every $n\geq 1$ and set $\eta_n = p_n^*\mu$, where we recall that $\mu$ is the $G$-invariant probability measure on $G/\Lambda$ coming from a Haar measure on $G$. Then $\eta_n$ is a probability measure on $G/\Lambda\times G/\Lambda$ which satisfies the three properties of Proposition \ref{propIoa} ($\eta_n$ is even $\Gamma$-invariant). Yet, $\eta_n(\Delta_X)=0$ for all $n$. Hence, the action $[\Gamma_1,\Gamma_1]\curvearrowright (X,\mu)$ does not have the property (T) relative to the space.

Since $[\Gamma_1,\Gamma_1]$ is a co-amenable subgroup the action of $\Gamma$ on $(X,\mu)$ does not have the property (T) relative to the space.

\medskip

In any case, the action $\Gamma\curvearrowright (X,\mu)$ does not have the property (T) relative to the space. 

\begin{flushright}
$\square$
\end{flushright}

\begin{rk}
The proof of implication $(ii)\Rightarrow (i)$ of Theorem \ref{thmp} shows that this implication holds when $\mu$ is replaced by any $\Gamma$-invariant atomless probability measure on $G/\Lambda$. 
\label{rkthmp}
\end{rk}

\section{Construction of rigid actions for linear groups}

\subsection{Preliminaries}

The following lemma will allow us to show that certain actions are free. 

\begin{lem}
Let $k$ be a local field of characteristic zero. Let $G$ be a Zariski-connected algebraic group over $k$. Let $\Lambda$ be a lattice in $G$. Let $\mu_{G/\Lambda}$ be the $G$-invariant probability measure on $G/\Lambda$. Let $g\in G$ be such that$$\mu_{G/\Lambda}(\{\overline{h}\in G/\Lambda\,:\,g\overline{h}=\overline{h}\})>0$$ then $g$ is contained in the centre $Z(G)$ of $G$.
\label{lemfree}
\end{lem}

\textbf{Proof:}

Let $X$ be a Borel fundamental domain for the action of $\Lambda$ on $G$ by right translations. Let $\mu$ be the $G$-invariant measure on $X$ associated to the Haar measure on $G$; we can assume that $\mu$ is a probability measure.

It is enough to show that $\mu(\{h\in X\,:\, h^{-1}gh\in \Lambda\})>0\Rightarrow g\in Z(G)$. 

Since $\Lambda$ is countable, there exists $\lambda\in\Lambda$ such that $$\mu(\{h\in X\,:\,h^{-1}gh=\lambda\})>0.$$ 
Set $A:=\{h\in X\,:\, h^{-1}gh=\lambda\})\subset G$. Since $A$ is a measurable subset of $G$ of positive measure, it is well known that $AA^{-1}$ contains a neighborhood of the identity in $G$.

\medskip

Observe that $AA^{-1}$ is contained in the stabiliser Stab$_G(g)$ of $g$ in $G$. Indeed, for $a,b\in A$ we have: $$ab^{-1}gba^{-1}=a\lambda a^{-1}=g.$$ Therefore, Stab$_G(g)$ contains an open non-empty subset. On the other hand, Stab$_G(g)$ is Zariski-closed. Hence Stab$_G(g)=G$ and this means that $g\in Z(G)$.

\begin{flushright}
$\square$
\end{flushright}

The proof of Proposition \ref{lemdecom} is based on the following result.

\begin{lem}
Let $\Gamma\subset$ GL$_n(k)$ be a finitely generated group over a local field $k$ of characteristic zero. Assume that $\Gamma$ has trivial solvable radical and that the Zariski-closure $G$ of $\Gamma$ in GL$_n(k)$ is simple. Then there exists a local field $\tilde{k}$ of characteristic zero and an injective morphism $$j:\Gamma\rightarrow \text{ GL}_n(\tilde{k})$$with unbounded image such that the Zariski-closure $S$ of $j(\Gamma)$ is simple.
\end{lem}

\textbf{Proof:}

First note that upon replacing $G$ by $G/Z(G)$, we can assume that $G$ is centre-free since $\Gamma$ has trivial solvable radical and so $\Gamma$ embeds in the centre-free linear group $G/Z(G)$ over $k$. As such, we can also assume that $G$ is simple as an abstract group (see chapter 29.5 in \cite{Hum}).

By Tits' Lemma (Lemma 4.1 \cite{Tits}), there exists a local field $\tilde{k}$ of characteristic zero and an injective morphism $$i:\Gamma\rightarrow \text{ GL}_n(\tilde{k})$$ with unbounded image.

Let $\tilde{G}$ be the Zariski-closure of $i(\Gamma)$ in GL$_n(\tilde{k})$ and let Rad$(\tilde{G})$ be the solvable radical of $\tilde{G}$. The projection $\Gamma\rightarrow \tilde{G}/$Rad$(\tilde{G})$ is injective since $\Gamma\cap$Rad$(\tilde{G})$ is trivial. So, as above, we can assume that $\tilde{G}$ is semi-simple.

We can write $\tilde{G}=S_{nc,1}\times\dots\times S_{nc,r}\times S_{c,1}\times S_{c,q}$ since $\Gamma$ is centre-free, where $S_{nc,l}$, $l=1,\dots,r$ are simple non-compact groups and $S_{c,l}$, $l=1,\dots,q$ are simple compact groups. As $i(\Gamma)$ is unbounded, there exists $l\in \{1,\dots,r\}$ such that the projection of $i(\Gamma)$ in $S_{nc,l}$ is unbounded. We may assume that $l=1$. 

We claim that the projection $$\pi_1:\Gamma\rightarrow S_{nc,1}$$ is injective.

Assume by contradiction that Ker$\pi_1=\Gamma_0\subset \Gamma$ is not trivial. Let $H$ be the Zariski closure of $\Gamma_0$ in GL$_n(k)$. Then $H$ is a non-trivial normal subgroup of $G$ and hence $H=G$.

Let $k_0$ be the field generated by the coefficients of the matrices representing a generating set of $\Gamma$ in GL$_n(k)$. Since $k$ is a field extension of $k_0$, $\Gamma_0$ and $\Gamma$ have the same Zariski-closure in GL$_n(k_0)$; $\overline{\Gamma_0}^{Zar,k_0}=\overline{\Gamma}^{Zar,k_0}$.

On the oher hand, $\tilde{k}$ being a field extension of $k_0$, the Zariski closures of $\Gamma$ and $\Gamma_0$ are also equal; $$\overline{\Gamma_0}^{Zar,\tilde{k}}=\overline{\Gamma}^{Zar,\tilde{k}}=\tilde{G}=S_{nc,1}\times\dots\times S_{nc,r}\times S_{c,1}\times S_{c,q}.$$ But $\Gamma_0\subset S_{nc,2}\times\dots\times S_{nc,r}\times S_{c,1}\times S_{c,q}$, so the Zariski-closures of $\Gamma_0$ and $\Gamma$ must be distinct. We get a contradiction and $\Gamma_0=\{e\}$. 

Therefore we obtain an injective morphism $$j:\Gamma\rightarrow S_{nc,1}$$with unbounded Zariski-dense image.

\begin{flushright}
$\square$
\end{flushright}

\textbf{Proof of Proposition \ref{lemdecom}:}

Let $G\subset$ GL$_n(k)$ be the Zariski-closure of $\Gamma$ in GL$_n(k)$. Again, we can assume that $G$ is semi-simple. 

We can write $G=\tilde{S}_1\times\dots\times\tilde{S}_r $, where $\tilde{S}_i$ is simple for $i=1,\dots,r$. For every $i=1,\dots,r$, the projection $\Gamma_i$ of $\Gamma$ in $\tilde{S}_i$ satisfies the conditions of the previous lemma.

Therefore, for $i=1,\dots,r$, we get a local field $k_i$ of characterisitic zero and an injective morphism $j_i:\Gamma_i\rightarrow$ GL$_n(k_i)$ with unbounded image and simple Zariski-closure $S_i$. 

Hence, the following morphism $$j:\Gamma\rightarrow S_1\times \dots S_r$$is injective and  has unbounded image and the projection of $j(\Gamma)$ on each factor $S_i$ is Zariski-dense. 

\begin{flushright}
$\square$
\end{flushright}

\subsection{Proof of Theorems \ref{thmr} and \ref{thml}}

\textbf{Proof of Theorem \ref{thmr}}

Let $\Gamma\subset$ GL$_n(k)$ be a non-amenable linear subgroup over a local field of characteristic zero $k$ with Rad$(\Gamma)=\{e\}$. By Proposition \ref{lemdecom}, there exist local fields $k_1,\dots,k_r$ of characteristic zero, and simple Lie groups $S_1,\dots,S_r$ over $k_1,\dots,k_r$ respectively such that $$\Gamma\hookrightarrow S_1\times\dots\times S_r$$with unbounded Zakiski-dense image on each factor.

For each factor $S_i$, $i=1,\dots,r$, let $\Lambda_i$ be a lattice (\cite{Bor}) and let $\mu_i$ be the $S_i$-invariant probability measure on $S_i/\Lambda$. Consider the natural action of $\Gamma$ on $(X,\mu)=(\prod_{i=1}^r S_i/\Lambda_i,\otimes_{i=1}^r \mu_i)$. Thanks to Proposition \ref{lemcoind} we can assume that the groups $S_i$ are Zariski-connected since the Zariski-component of the identity is a normal finite index subgroup. By Lemma \ref{lemfree} this action is free. 

In order to show that this action is ergodic and has the property (T) relative to the space, it is enough to show that for $i=1,\dots,r$ the action $\Gamma_i\curvearrowright(S_i/\Lambda)$ is ergodic and has the property (T) relative to the space, according to Lemma \ref{lemPr}.

By contradiction, assume that this action does not have the property (T) relative to the space. Then, by Theorem \ref{thmp}, there exists a $\Gamma$-invariant probability measure on $\mathbb{P}(\frak{s_i})$ where $\frak{s_i}$ is the Lie algebra of $S_i$. 

The subgroup Ad$(\Gamma)\subset$ GL$(\frak{s_i})$ is unbounded since Ad$:S\rightarrow $Aut$(\frak{s_i})^0$ is an isomorphism with non-compact image (see e.g. chapter 14.1 in \cite{Hum}). Hence, by Furstenberg's Lemma, there exists a finite index subgroup $\Gamma_0\subset\Gamma$ which leaves a proper subspace $V_0\subset\frak{s_i}$ invariant. Since $\Gamma$ is Zariski-dense in $S_i$, $\Gamma_0$ is also Zariski-dense in $S_i$.

As such, $V_0$ is Ad$(S_i)$-invariant. Hence, $V_0$ is a proper ideal of $\frak{s_i}$. This contradicts the simplicity of $S_i$.

Therefore the action $\Gamma\curvearrowright (S_i/\Lambda,\mu)$ has the property (T) relative to the space. 

Regarding ergodicity, it is an application of Howe-Moore theorem (\cite{HM}). 

It ends the proof.

\begin{flushright}
$\square$
\end{flushright}

\textbf{Proof of Theorem \ref{thml}}

It is sufficient to consider the projection $$\pi:\Gamma\rightarrow \Gamma/\text{Rad}(\Gamma).$$ Since $\pi(\Gamma)$ satisfies the assumptions of Theorem \ref{thmr}, $\Gamma$ clearly has an ergodic action with the property (T) relative to the space.

\begin{flushright}
$\square$
\end{flushright}

\subsection{Application}

In \cite{Ioa}, Ioana and Shalom asked whether $\mathbb{F}_2\times\Z$ admits a free ergodic rigid action. We claim that it is indeed the case. More generally we prove the following result.

\begin{prop}
Let $\Gamma$ be a finitely-generated non-amenable linear group which has trivial solvable radical. 
Assume that there exists a morphism $\varphi:\Gamma\rightarrow$ SL$_n(\R)$ such that the $\R$-Zariski-closure $\overline{\varphi(\Gamma)}^Z(\R)$ is non-amenable.

Then $\Gamma\times \Z$ admits a free ergodic action with the property (T) relative to the space.
\label{propmoi}
\end{prop}




For the proof, we will use the following result of T.Fernos.

\begin{thm}[\cite{Fer}]
Let $\Gamma$ be a finitely generated group. The following assertions are equivalent:
\begin{itemize}
\item[(i)] there exists a morphism $\varphi:\Gamma\rightarrow$ SL$_n(\R)$ such that the $\R$-Zariski-closure $\overline{\varphi(\Gamma)}^Z(\R)$ is non-amenable,
\item[(ii)] there exists an abelian group $A$ of nonzero finite $\Q$-rank and a morphism $\varphi' : \Gamma\rightarrow$ Aut$(A)$ such that the corresponding group pair $(A\rtimes_{\varphi'}\Gamma,A)$ has the relative property (T). 
\end{itemize}
\end{thm}

\begin{rk}\begin{itemize}\item[(i)] The abelian group $A$ appearing in (ii) above, can be taken of the form $\Z[S^{-1}]^N$ where $S$ is some finite set of rational primes.
\item[(ii)] The action $\Gamma\curvearrowright (\hat{A},\mu)$ has the property (T) relative to the space (see Proposition 5.1 \cite{Pop})
\end{itemize}
\end{rk}

\textbf{Proof of Proposition \ref{propmoi}}

By Fernos' Theorem, there exists a morphism $\varphi': \Gamma\rightarrow $ SL$_n(\Z[S^{-1}])$ such that the pair $(A\rtimes \varphi'(\Gamma),A)$ has relative property (T), with $A=\Z[S^{-1}]^n$. So, the associated action $\Gamma\curvearrowright(\hat{A},\lambda)$ has the property (T) relative to the space, where $\lambda$ is the normalized Haar measure on the compact group $\hat{A}$. By Remark \ref{rkthmp}, upon replacing $\mu$ by some of its ergodic components, we may assume that $\mu$ is ergodic for the $\Gamma$-action. 

\medskip
 By Lemma \ref{lemPr}, the diagonal action of $\Gamma$ on $(\hat{A}\times \hat{A},\lambda\otimes\lambda)$ has the property (T) relative to the space.

On the other hand, $\Z$ acts freely on $(\hat{A}\times\hat{A},\lambda\otimes\lambda)$ by: $$m.(a,b)=(a+mb,b),\quad \text{for all }m\in \Z, (a,b)\in \hat{A}\times\hat{A}$$where the group law on $\hat{A}$ is written additively. The action $\Gamma\times \Z\curvearrowright(\hat{A}\times\hat{A},\lambda\otimes\lambda),$ given by $$ (\gamma,m).(a,b)=(\gamma(a+mb),\gamma(b)),$$ has the property (T) relative to the space.

Now by Theorem $\ref{thmr}$, $\Gamma$ admits a free ergodic action $\Gamma\curvearrowright(X,\mu)$ which has the property (T) relative to the space. 

Finally the action $\Gamma\times\Z \curvearrowright (X\times \hat{A}\times\hat{A},\mu\otimes\lambda\otimes\lambda)$ given for all $\gamma\in\Gamma, m\in\Z, x\in X, (a,b)\in \hat{A}\times\hat{A}$ by $$(\gamma,m).(x,a,b)=(\gamma.x,\gamma(a+mb),\gamma(b))$$ is free ergodic and has the property (T) relative to the space.

\begin{flushright}
$\square$
\end{flushright}

\begin{rk}
The proof of Proposition \ref{propmoi} shows that the following more general result is true. Let $\Gamma_1$ be a subgroup of SL$_n(\Z[S^{-1}])$ such that the action $\Gamma_1\curvearrowright ((\R\times \prod_{p\in S} \Q_p)^n / \Z[S^{-1}]^n ,\lambda)$ is free ergodic and has the property (T) relative to the space, where $S$ is a finite subset of rational primes. Then, for any subgroup $\Gamma_2 \subset$ SL$_m(\Z[S^{-1}])$ the group $\Gamma_1\times \Gamma_2$ admits a free ergodic action with the property (T) relative to the space.
\end{rk}



\bibliographystyle{alpha}

\begin{small}
\bibliography{biblio}
\end{small}

\Addresses

\end{document}